\documentclass[a4paper,reqno]{amsart} 
\usepackage[T1]{fontenc}
\usepackage{lmodern}
\usepackage{mathrsfs} 
\usepackage{amsmath,amssymb,amsthm} 
\usepackage{mathtools}

\usepackage{upref}
\usepackage{enumerate} 
\usepackage{indentfirst}
\usepackage[numbers,sort&compress]{natbib}

\usepackage{hyperref}%
\hypersetup{
hypertexnames=false, 
linkbordercolor=[rgb]{0.5, 1, 1},%
urlbordercolor=[rgb]{0.5, 1, 1},%
citebordercolor=[rgb]{0.5, 1, 0.5}}%
\usepackage{graphicx} 
\usepackage{subcaption} 

\newtheorem{theorem}{Theorem} 
\newtheorem{lemma}{Lemma} 
\newtheorem{proposition}{Proposition} 

\theoremstyle{definition} 
\newtheorem{problem}{Problem} 

\theoremstyle{remark}
\newtheorem{remark}{Remark} 

\newcommand\setA{\mathscr{A}}
\newcommand\setB{\mathscr{B}}
\newcommand\setM{\mathscr{M}}
\newcommand\bbN{\mathbb{N}}
\begin{document}

\title[The problem of boundedness of matrix products]{The problem of
  boundedness of matrix products with alternating factors from two sets of
  matrices revisited}

\author{Victor Kozyakin}

\thanks{The research is supported by MSHE RF GZ project FSMG-2026-0020}

\address{Higher School of Modern Mathematics MIPT,
  9 Institutskiy per., Dolgoprudny, Moscow Region, 141701, Russian Federation}

\email{koziakin.vs@mipt.ru}


\date{}

\begin{abstract}\sloppy
  The paper addresses the boundedness of matrix products $A_{n}B_{n}\cdots
    A_{1}B_{1}$---where the factors are drawn from two sets of matrices,
  $A_{i}\in\setA$ and $B_{i}\in\setB$---achieved through a suitable choice
  of the matrices $\{B_{i}\}$. It is assumed that for any sequence of
  matrices $\{A_{i}\}$, there exists a sequence of matrices $\{B_{i}\}$ such
  that the sequence of matrix products $\{A_{n}B_{n}\cdots
    A_{1}B_{1}\}_{n=1}^{\infty}$ is bounded in norm. It is proved that this
  boundedness holds, without additional conditions, for the class of
  so-called prefix sequences $\{A_{i}\}$ characterized by a specific fixed
  ``prefix''---an initial segment of the sequence $\{A_{i}\}$.
\end{abstract}

\maketitle

\section{Introduction}\label{S:intro}

Let $\setM$ be a set of matrices. We denote by $\setM^{\mathbb{N}}$ the
set of infinite sequences of matrices $\{M_{i}\in\setM\colon
  i\in\mathbb{N}\}$, and by $\setM(p,q)$ the space of real $(p\times
  q)$-matrices equipped with the topology of entrywise convergence.

Let $\setA\subset\setM(N,N)$ be a finite set of (real or complex)
matrices. As shown in~\cite{Koz:AiT90:10:e} (and somewhat later, in
various modified forms, in~\cite{DaubLag:LAA92, BerWang:LAA92,
  Gurv:LAA95, Hartfiel:02}), if for every sequence of matrices
$\{A_{n}\}\in\setA^{\bbN}$ the sequence of norms $\{\|A_{n}\cdots
  A_{1}\|\}_{n=1}^{\infty}$ converges to zero, then this convergence to
zero is, in fact, uniform with respect to the choice of the matrix
sequence. Similarly, if for every sequence of matrices
$\{A_{n}\}\in\setA^{\bbN}$ the sequence of norms $\{\|A_{n}\cdots
  A_{1}\|\}_{n=1}^{\infty}$ is bounded, then this boundedness is also
uniform with respect to the choice of the matrix sequence. Precise
statements of these results can be found, for example,
in~\cite{Koz:AiT90:10:e}.

In~\cite{Koz:DDNS18}, an attempt was made to extend the convergence
analysis of matrix products to products involving alternating factors
from two sets of matrices. The need for such an analysis was motivated
by the problem of stabilizability for so-called discrete-time switched
linear systems; see, for example,~\cite{Stanford:SIAMJCO79,
  Koz:AiT90:10:e, SunGe05, JunMas:SIAMJCO17, LinAnt:IEEETAC09}. Among
these publications, we would like to highlight the
work~\cite{JunMas:SIAMJCO17}, in which alternative approaches are
considered and discussed.

More precisely, the following question was considered
in~\cite{Koz:DDNS18}: given finite sets of matrices
$\setA\subset\setM(N, M)$ and $\setB\subset \setM(M, N)$, what can be
said about the rate of convergence to zero of the matrix products
\[
  A_{n}B_{n}\cdots A_{1}B_{1},\qquad A_{i}\in\setA,~B_{i}\in\setB,
\]
under the assumption that the convergence $A_{n}B_{n}\cdots
  A_{1}B_{1}\to 0$ can be ensured for any sequence of matrices
$\{A_{i}\}$ by a suitable choice of the sequence of matrices
$\{B_{i}\}$? As it turned out~\cite[Theorem~2]{Koz:DDNS18}, the
following assertion holds true, which naturally generalizes similar
results from~\cite{DaubLag:LAA92, BerWang:LAA92, Gurv:LAA95,
  Hartfiel:02} concerning the convergence of matrix products with factors
drawn from a single set of matrices.

\begin{proposition}\label{St:1}
  Let $\setA\subset\setM(N,M)$ and $\setB\subset\setM(M,N)$ be finite sets of matrices with the property that
  \[
    \forall \{A_{n}\}\in\setA^{\bbN}\quad\exists \{B_{n}\}\in\setB^{\bbN}:\quad
    A_{n}B_{n}\cdots A_{1}B_{1}\to0\quad\text{as}\quad n\to\infty.
  \]
  Then
  \[
    \exists C<\infty:\quad\forall \{A_{n}\}\in\setA^{\bbN}\quad
    \exists \{B_{n}\}\in\setB^{\bbN}:\quad\|A_{n}B_{n}\cdots A_{1}B_{1}\|\le
    C\lambda^{n}\quad\forall n\in\bbN.
  \]
\end{proposition}

The question of the boundedness of matrix products with factors drawn
from two sets of matrices turned out to be more complex. A partial
answer to this question is provided in~\cite[Thms.~1,~2]{Koz:ArXiv20}:

\begin{proposition}\label{St:2}
  Let $\setA\subset\setM(N,M)$ and $\setB\subset\setM(M,N)$ be finite sets of matrices with the property that
  \begin{equation}\label{eq:assumption}
    \forall \{A_{n}\}\in\setA^{\bbN}\quad\exists \{B_{n}\}\in\setB^{\bbN}:\quad
    \sup_{n\in\mathbb{N}}\|A_{n}B_{n}\cdots A_{1}B_{1}\|<\infty.
  \end{equation}
  Suppose, furthermore, that one of the following conditions holds:
  \begin{itemize}
    \item[$(i)$] the pairs of matrix sets $\setA$ and $\setB$ are
      non-singular:
      \begin{align*}
        \det(AB)\neq 0,\quad \forall
        A\in\setA,~B\in\setB,\qquad\text{if}~N\le M, \\
        \det(BA)\neq 0,\quad \forall  A\in\setA,~B\in\setB,\qquad\text{if}~N> M,
      \end{align*}

    \item[$(ii)$] the pairs of matrix sets $\setA$ and $\setB$ are
      non-negative:
      \[
        a_{ij}\ge0,\quad b_{ij}\ge0,\qquad \forall
        A=(a_{ij})\in\setA,~B=(b_{ij})\in\setB,
      \]
      and, furthermore, each matrix in $\setA$ and $\setB$ is
      row-non-singular---i.e., it cannot have rows consisting entirely of
      zero elements:
      \[
        \min_{i}\sum_{j}a_{ij}>0,\quad\min_{i}\sum_{j}b_{ij}>0,\qquad
        \forall A=(a_{ij})\in\setA,~B=(b_{ij})\in\setB.
      \]
  \end{itemize}

  Then
  \begin{equation}\label{eq:mainineq}
    \exists C<\infty:\quad\forall \{A_{n}\}\in\setA^{\bbN}\quad
    \exists \{B_{n}\}\in\setB^{\bbN}:\quad
    \sup_{n\in\mathbb{N}}\|A_{n}B_{n}\cdots A_{1}B_{1}\|\le C.
  \end{equation}
\end{proposition}

For brevity, we shall henceforth refer to statements of the
type~\eqref{eq:mainineq} as the \emph{uniform boundedness principle for
  matrix products with alternating factors}.

The goal of this paper is to obtain two additional conditions under
which the aforementioned principle of uniform boundedness holds,
thereby supplementing Prorosition~\ref{St:2}. In the first of
these---Theorem~\ref{thm:main}---it is proven that, under
assumption~\eqref{eq:assumption}, the uniform boundedness of matrix
products holds (without any further conditions) for the class of
so-called ``prefix'' sequences $\{A_{i}\}$, which are defined by a
fixed ``prefix''---an initial segment of the sequence $\{A_{i}\}$. The
second result, Theorem~\ref{thm:tailwise}, demonstrates that the
situation simplifies significantly if, rather than requiring the
boundedness of initial products of the form $A_{n}B_{n} \cdots
  A_{1}B_{1}$, one instead requires the uniform boundedness---with
respect to the starting index $m$---of the so-called ``tail'' products
$A_{n}B_{n}\cdots A_{m}B_{m}$.

\section{Main result}\label{sec:main}

Let us recall the definitions and results needed later. Let a set of
matrices $\setA \subset \setM(N, M)$ and a finite ordered collection of
matrices
\begin{equation}\label{eq:1}
  \boldsymbol{X}=\{X_{1},\ldots,X_{k}\},\quad X_{i}\in\setA
\end{equation}
be given. We denote by $\setA^{\bbN}_{\boldsymbol{X}}$ the set of all
\emph{$\boldsymbol{X}$-prefix} sequences $\{A_{n}\} \in \setA^{\bbN}$, i.e.,
the collection of sequences $\{A_{n}\} \in \setA^{\bbN}$ satisfying the
condition
\[
  A_{i}=X_{i},\quad i=1,2,\ldots,k.
\]
\begin{theorem}\label{thm:main}
  Let $\setA\subset\setM(N,M)$ and $\setB\subset\setM(M,N)$ be finite sets of 
  matrices for which condition~\eqref{eq:assumption} holds. 

  Then, for a certain collection of matrices~\eqref{eq:1}, the principle
  of uniform boundedness for matrix products with alternating factors
  $\{A_{n}\}\in\setA^{\bbN}_{\boldsymbol{X}}$ and $\{B_{n}\}\in\setB$
  holds; that is,
  \begin{equation}\label{eq:tailineq}
    \exists C<\infty:\quad\forall \{A_{n}\}\in\setA^{\bbN}_{\boldsymbol{X}}\quad
    \exists \{B_{n}\}\in\setB^{\bbN}:\quad
    \sup_{n\in\mathbb{N}}\|A_{n}B_{n}\cdots A_{1}B_{1}\|\le C.
  \end{equation}
\end{theorem}

Before proceeding directly to the proof of this theorem, let us recall
the topological properties of sets of sequences and establish an
auxiliary statement.

For any finite set of matrices $\setM$, we define a (discrete) metric
$d(\cdot,\cdot)$ on it by setting, for any $M, \tilde{M} \in \setM$,
\[
  d(M,\tilde{M})=
  \begin{cases} 0, & \text{if } M=\tilde{M},    \\
                1, & \text{if } M\neq\tilde{M}.
  \end{cases}
\]
We also define a metric $\rho(\cdot,\cdot)$ on the space of sequences
$\setM^{\bbN}$ by setting, for any pair of sequences
$\{M_{i}\}_{i=1}^{\infty},
  \{\tilde{M}_{i}\}_{i=1}^{\infty}\in\setM^{\bbN}$,
\[
  \rho(\{M_{i}\},\{\tilde{M}_{i}\})=\sum_{n=1}^{\infty}2^{-n}d(M_{n},\tilde{M}_{n}).
\]
Since $d(M,\tilde{M})\le1$ for any matrices $M,\tilde{M}\in\setM$, the
metric $\rho(\cdot,\cdot)$ is well-defined, and for any elements
(sequences of matrices) in $\setM^{\bbN}$, it also does not exceed $1$.

The space $\setM^{\bbN}$ can be viewed as the direct (Cartesian)
product of a countable number of discrete spaces $\setM$, endowed with
the Tychonoff product topology. In fact, the space $\setM^{\bbN}$ is
metrizable---the metric $\rho(\cdot,\cdot)$ on it induces the same
Tychonoff product topology. Convergence in this topology is simply
coordinate-wise convergence of sequences of elements of $\setM$ with
respect to the metric $d(\cdot,\cdot)$, which---given the finiteness of
the set $\setM$---amounts to the eventual stabilization of the
corresponding coordinates. Moreover, for finite sets $\setM$, the space
$\setM^{\bbN}$---equipped with the metric $\rho(\cdot,\cdot)$ (or,
equivalently, the product topology)---is complete and compact by
Tychonoff's theorem.

Formally, we defined the metrics $d(\cdot,\cdot)$ and
$\rho(\cdot,\cdot)$ above for the sets $\setM$ and $\setM^{\bbN}$,
respectively. Subsequently, when we need to consider other sets of
matrices ($\setA$, $\setB$, etc.), we will not change the notation for
these metrics, provided this does not lead to ambiguity.

For a set of matrices $\setA$ and an arbitrary constant $C>0$, we
define the set
\[
  \setA^{\bbN}(C)\subseteq\setA^{\bbN}
\]
as the collection of those sequences
$\{A_{i}\}_{i=1}^{\infty}\in\setA^{\bbN}$ for which there exist
sequences $\{B_{i}\}_{i=1}^{\infty}\in\setB^{\bbN}$ such that
\[
  \|A_{i}B_{i}\cdots A_{1}B_{1}\|\le C,\qquad i\in\bbN.
\]
\begin{lemma}\label{lem:1}
  For every $C>0$, the set $\setA^{\bbN}(C)$ is closed in the topology of the
  space~$\setA^{\bbN}$.
\end{lemma}

\begin{proof}
  Only the case where the set $\setA^{\bbN}(C)$ is non-empty requires proof.
  Let $\boldsymbol{A^{*}}=\{A^{*}_{n}\}_{n=1}^{\infty}\in\setA^{\bbN}$ be an
  element of the closure of the set $\setA^{\bbN}(C)$; that is, for some
  sequence
  \begin{equation}\label{eq:Ak}
    \boldsymbol{A^{(k)}}=\{A^{(k)}_{n}\}_{n=1}^{\infty}\in\setA^{\bbN}(C),\qquad
    k\in\bbN,
  \end{equation}
  the equality
  \[
    \boldsymbol{A^{*}}=\lim_{k\to\infty}\boldsymbol{A^{(k)}}
  \]
  holds, where the limit is understood in the sense of convergence in the
  product topology on $\setA^{\bbN}$.

  By the definition of the set $\setA^{\bbN}(C)$, relation~\eqref{eq:Ak}
  implies the existence of sequences
  \[
    \boldsymbol{B^{(k)}}=\{B^{(k)}_{n}\}_{n=1}^{\infty}\in\setB^{\bbN},\qquad
    k\in\bbN,
  \]
  such that for each $k\in\bbN$, the inequalities
  \begin{equation}\label{eq:ABprod}
    \|A^{(k)}_{i}B^{(k)}_{i}\cdots A^{(k)}_{1}B^{(k)}_{1}\|\le C,\qquad i\in\bbN,
  \end{equation}
  hold.

  Let us now note that the space $\setB^{\bbN}$ is compact, as previously
  mentioned. Therefore, a subsequence of the sequence
  $\{\boldsymbol{B^{(k)}}\}$ converges to some limit. Without loss of
  generality, we may assume that the sequence $\boldsymbol{B^{(k)}}$
  itself converges to a limit
  $\boldsymbol{B^{*}}=\{B^{*}_{n}\}_{n=1}^{\infty}\in\setB^{\bbN}$, i.e., \[
    \boldsymbol{B^{*}}=\lim_{k\to\infty}\boldsymbol{B^{(k)}},
  \]
  where the limit is understood in the sense of convergence in the
  product topology on $\setB^{\bbN}$.

  Since convergence in the product topology of the spaces $\setA^{\bbN}$
  and $\setB^{\bbN}$ is equivalent to coordinate-wise convergence of
  sequences (as noted earlier), we can pass to the limit
  in~\eqref{eq:ABprod}:
  \[
    \|A^{*}_{i}B^{*}_{i}\cdots A^{*}_{1}B^{*}_{1}\|=
    \lim_{k\to\infty}\|A^{(k)}_{i}B^{(k)}_{i}\cdots A^{(k)}_{1}B^{(k)}_{1}\|\le
    C,\qquad i\in\bbN.
  \]
  The resulting estimates imply that
  \[
    \boldsymbol{A^{*}}=\{A^{*}_{n}\}_{n=1}^{\infty}\in\setA^{\bbN}(C).
  \]

  The lemma is proved.
\end{proof}

We can now proceed to the proof of Theorem~\ref{thm:main}.

\begin{proof}[Proof of Theorem~\ref{thm:main}]
  By the hypothesis of the theorem, every sequence
  $\{A_{i}\}_{i=1}^{\infty}\in\setA^{\bbN}$ belongs to some set
  $\setA^{\bbN}(n)$, $n\in\mathbb{N}$. Therefore,
  \[
    \setA^{\bbN}=\bigcup_{n\ge1}\setA^{\bbN}(n).
  \]

  As noted earlier, the space $\setA^{\bbN}$ is complete in the topology
  of coordinate-wise convergence of sequences $\{A_{i}\}_{i=1}^{\infty}
    \in \setA^{\bbN}$, which is equivalent to its completeness in the
  metric $\rho(\cdot, \cdot)$. Then, by the Baire category theorem---and
  given that the sets $\setA^{\bbN}(n)$ are closed by
  Lemma~\ref{lem:1}---there exists an $n_{*} \in \mathbb{N}$ such that
  the set $\setA^{\bbN}(n_{*})$ contains interior points. In particular,
  there exists a sequence
  \[
    \{X_{i}\}_{i=1}^{\infty} \in \setA^{\bbN}(n_{*})
  \]
  that belongs to the set $\setA^{\bbN}(n_{*})$ together with some
  neighborhood of radius $\varepsilon$ in the metric $\rho(\cdot,
    \cdot)$.

  Let us now choose an arbitrary $k\in\mathbb{N}$ satisfying the relation
  \[
    2^{-k}\le\varepsilon,
  \]
  and estimate the distance in the metric $\rho(\cdot,\cdot)$ between an
  arbitrary element
  \begin{equation}\label{eq:defA}
    \{A_{i}\}\in\setA^{\bbN}_{\boldsymbol{X}},\quad\text{where}\quad
    \boldsymbol{X}=\{X_{1},\ldots,X_{k}\},
  \end{equation}
  and the sequence $\{X_{i}\}_{i=1}^{\infty}$:
  \[
    \rho(\{A_{i}\},\{X_{i}\})=\sum_{n=1}^{\infty}2^{-n}d(A_{n},X_{n})=
    \sum_{n=1}^{k}2^{-n}d(A_{n},X_{n})+\sum_{n=k+1}^{\infty}2^{-n}d(A_{n},X_{n}).
  \]

  Here, all terms in the first sum on the right-hand side are zero by
  virtue of~\eqref{eq:defA}, while all quantities $d(A_{n},X_{n})$ in the
  second sum on the right-hand side do not exceed $1$ by the definition
  of the metric $d(\cdot,\cdot)$. Therefore,
  \[
    \rho(\{A_{i}\},\{X_{i}\})\le\sum_{n=k+1}^{\infty}2^{-n}=2^{-k}\le\varepsilon.
  \]
  It follows from the obtained estimate that every
  sequence~\eqref{eq:defA} belongs to the $\varepsilon$-neighborhood of
  the sequence $\{X_{i}\}_{i=1}^{\infty}$, and thus belongs to the set
  $\setA^{\bbN}(n_{*})$ as well. Then, by the definition of the set
  $\setA^{\bbN}(\cdot)$, for any sequence of matrices
  $\{A_{n}\}\in\setA^{\bbN}_{\boldsymbol{X}}$, there exists a sequence of
  matrices $\{B_{n}\}\in\setB^{\bbN}$ such that $\|A_{n}B_{n}\cdots
    A_{1}B_{1}\|\le C:=n_{*}$ for all $n\in\bbN$.

  Theorem~\ref{thm:main} is proved.
\end{proof}

\begin{remark}\label{rem:1}
  Proposition~\ref{St:2} follows from Theorem~\ref{thm:main} when condition (i)
  holds. Indeed, under the conditions of Theorem~\ref{thm:main}, there are
  only a finite number of distinct prefix parts of the sequences
  $\{A_{n}\}\in\setA^{\bbN}_{\boldsymbol{X}}$,
  and when condition (i) of Proposition~\ref{St:2} is satisfied, the
  alternating products of matrices $A_{i}\in\setA$ and $B_{i}\in\setB$
  corresponding to these prefix parts are invertible. Then the matrix norms
  \[
    \|A_{n_{*}+k}B_{n_{*}+k}\cdots A_{n_{*}+1}B_{n_{*}+1}\|
  \]
  are bounded above, uniformly in $k\in\mathbb{N}$, by the norms of the
  matrices
  \[
    \|A_{n_{*}+k}B_{n_{*}+k}\cdots A_{1}B_{1}\|
  \]
  and are therefore uniformly bounded in total by Theorem~\ref{thm:main}.
  Since the set of all ``tail'' subsequences
  \[
    \{A_{i}\}_{i=n_{*}+1}^{\infty},\quad
    \text{where}\quad\{A_{i}\}\in\setA^{\bbN}_{\boldsymbol{X}},
  \]
  coincides with the entire space $\setA^{\bbN}$ by
  Theorem~\ref{thm:main}, this implies the validity of
  Proposition~\ref{St:2} under condition (i).
\end{remark}

\section{Uniform boundedness of the ``tails''}\label{sec:tailcase}

In Proposition~\ref{St:2} and Theorem~\ref{thm:main}, the fundamental
assumption~\eqref{eq:assumption} and the corresponding
statements~\eqref{eq:mainineq} and~\eqref{eq:tailineq} concerned the
boundedness of the norms of matrix products of the form
\[
  A_{n}B_{n}\cdots A_{1}B_{1},\qquad n=1,2,\ldots\,.
\]
The situation simplifies significantly if one requires the ``tail''
products
\[
  A_{n}B_{n}\cdots A_{m}B_{m},\qquad n=m,m+1,\ldots\,,
\]
to be uniformly bounded with respect to the starting index $m$.
\begin{theorem}\label{thm:tailwise}
  Let $\setA\subset\setM(N,M)$ and $\setB\subset\setM(M,N)$ be finite
  sets of matrices. Then the following two statements
  \begin{gather}\label{eq:assumption-new}
    \forall \{A_{n}\}\in\setA^{\bbN}\quad\exists \{B_{n}\}\in\setB^{\bbN}:\quad
    \sup_{n\ge m\ge1}\|A_{n}B_{n}\cdots A_{m}B_{m}\|<\infty,\\
    \shortintertext{and}
    \label{eq:tailineq-new}
    \exists C<\infty:\quad\forall \{A_{n}\}\in\setA^{\bbN}_{\boldsymbol{X}}\quad
    \exists \{B_{n}\}\in\setB^{\bbN}:\quad
    \sup_{n\ge m\ge1}\|A_{n}B_{n}\cdots A_{m}B_{m}\|\le C,
  \end{gather}
  are equivalent to each other.
\end{theorem}

The proof of this theorem follows the proof of Theorem~\ref{thm:main}
almost verbatim. Nevertheless, for the sake of completeness, we provide
the relevant details.

Given a set of matrices $\setA$ and an arbitrary constant $C>0$, we
define the set
\[
  \setA_{*}^{\bbN}(C)\subseteq\setA^{\bbN}
\]
as the collection of all sequences
$\{A_{i}\}_{i=1}^{\infty}\in\setA^{\bbN}$ for which there exist
sequences $\{B_{i}\}_{i=1}^{\infty}\in\setB^{\bbN}$ such that
\[
  \|A_{n}B_{n}\cdots A_{m}B_{m}\|\le C,\qquad \forall  n,m\in\mathbb{N}:~n\ge m\ge1.
\]
\begin{lemma}\label{lem:2}
  For every $C>0$, the set $\setA_{*}^{\bbN}(C)$ is closed in the topology of
  the space~$\setA^{\bbN}$.
\end{lemma}

The proof of this lemma follows the proof of Lemma~\ref{lem:1} almost
verbatim and is therefore omitted. Let us proceed to the proof of
Theorem~\ref{thm:tailwise}.

\begin{proof}[Proof of Theorem~\ref{thm:tailwise}]
  Statement~\eqref{eq:tailineq-new} of the theorem obviously implies statement~\eqref{eq:assumption-new}. Therefore, we shall prove only that statement~\eqref{eq:assumption-new} implies statement~\eqref{eq:tailineq-new}.

  According to the conditions of the theorem, every sequence
  $\{A_{i}\}_{i=1}^{\infty}\in\setA^{\bbN}$ belongs to some set
  $\setA_{*}^{\bbN}(n)$, where $n\in\mathbb{N}$. Consequently,
  \[
    \setA^{\bbN}=\bigcup_{n\ge1}\setA_{*}^{\bbN}(n).
  \]

  As noted earlier, the space $\setA^{\bbN}$ is complete under the
  topology of coordinate-wise convergence of sequences
  $\{A_{i}\}_{i=1}^{\infty}\in\setA^{\bbN}$, which is equivalent to
  saying that it is also complete with respect to the metric
  $\rho(\cdot,\cdot)$. By the Baire category theorem---and given that the
  sets $\setA_{*}^{\bbN}(n)$ are closed (as established in
  Lemma~\ref{lem:2})---there exists an $n_{*}\in\mathbb{N}$ such that the
  set $\setA_{*}^{\bbN}(n_{*})$ contains interior points. Specifically,
  this means there exists a sequence
  \[
    \{X_{i}\}_{i=1}^{\infty}\in\setA_{*}^{\bbN}(n_{*}),
  \]
  which belongs to the set $\setA_{*}^{\bbN}(n_{*})$ along with some
  neighborhood of radius $\varepsilon$ defined with respect to the metric
  $\rho(\cdot,\cdot)$.

  Let us now choose an arbitrary $k\in\mathbb{N}$ satisfying the relation
  \[
    2^{-k}\le\varepsilon,
  \]
  and estimate the distance, in terms of the metric $\rho(\cdot,\cdot)$,
  between an arbitrary element
  \begin{equation}\label{eq:defA0}
    \{A_{i}\}\in\setA^{\bbN}_{\boldsymbol{X}},\quad\text{where}\quad
    \boldsymbol{X}=\{X_{1},\ldots,X_{k}\},
  \end{equation}
  and the sequence $\{X_{i}\}_{i=1}^{\infty}$:
  \[
    \rho(\{A_{i}\},\{X_{i}\})=\sum_{n=1}^{\infty}2^{-n}d(A_{n},X_{n})=
    \sum_{n=1}^{k}2^{-n}d(A_{n},X_{n})+\sum_{n=k+1}^{\infty}2^{-n}d(A_{n},X_{n}).
  \]
  Here, all terms in the first summation on the right-hand side are equal
  to zero by virtue of~\eqref{eq:defA0}, while all quantities
  $d(A_{n},X_{n})$ in the second summation do not exceed $1$, by the
  definition of the metric $d(\cdot,\cdot)$. Therefore,
  \[
    \rho(\{A_{i}\},\{X_{i}\})\le\sum_{n=k+1}^{\infty}2^{-n}=2^{-k}\le\varepsilon.
  \]
  It follows from this estimate that every sequence of the
  form~\eqref{eq:defA0} belongs to the $\varepsilon$-neighborhood of the
  sequence $\{X_{i}\}_{i=1}^{\infty}$, and, consequently, also belongs to
  the set $\setA_{*}^{\bbN}(n_{*})$. Then, by the definition of the set
  $\setA_{*}^{\bbN}(\cdot)$, for any sequence of matrices
  $\{A_{n}\}\in\setA^{\bbN}_{\boldsymbol{X}}$, there exists a sequence of
  matrices $\{B_{n}\}\in\setB^{\bbN}$ such that $\|A_{n}B_{n}\cdots
    A_{m}B_{m}\|\le C:=n_{*}$ for all $n,m\in\mathbb{N}$ satisfying the
  inequalities $n\ge m\ge1$.

  It follows that for every sequence of matrices $\{A_n\} \in
    \mathcal{A}^{\mathbb{N}}_{\boldsymbol{X}}$, there exists a sequence of
  matrices $\{B_n\} \in \mathcal{B}^{\mathbb{N}}$ such that for the
  $n_*$-tails of these sequences, defined as
  \[
    \{\bar{A}_{i}\}_{i=1}^{{\infty}}=\{A_{i+n_{*}}\}_{i=1}^{{\infty}},\quad
    \{\bar{B}_{i}\}_{i=1}^{{\infty}}=\{B_{i+n_{*}}\}_{i=1}^{{\infty}}
  \]
  the following inequalities hold:
  \[
    \|\bar{A}_n \bar{B}_n \cdots \bar{A}_m \bar{B}_m\| \le C, \qquad \forall n, m
    \in \mathbb{N} \text{ such that } n \ge m \ge 1.
  \]
  It remains to note that the set of all $n_*$-tails $\{\bar{A}_n\}$ of
  matrix sequences $\{A_n\} \in
    \mathcal{A}^{\mathbb{N}}_{\boldsymbol{X}}$ coincides with the set of
  all sequences in $\mathcal{A}^{\mathbb{N}}$, which implies
  statement~\eqref{eq:tailineq-new} of the theorem.

  Theorem~\ref{thm:tailwise} is proved.
\end{proof}

\section{Open Problems}\label{sec:problems}

As proven in~\cite{Koz:AiT90:10:e}, if, for a given finite set of
matrices $\setA\subset\setM(N,N)$, every sequence of matrices
$\{A_{n}\}\in\setA^{\bbN}$ is bounded in norm, then this boundedness is
uniform with respect to the choice of the matrix sequence. It follows,
in particular, that in this case, all finite matrix products
${A_{n}\cdots A_{2}A_{1}}$ are also uniformly bounded in norm.

In this regard, and also in view of the results from
Sections~\ref{sec:main} and~\ref{sec:tailcase}, the following questions
arise:
\begin{problem}
  Let $\setA\subset\setM(N,M)$ and $\setB\subset\setM(M,N)$ be finite sets of
  matrices that satisfy assumption~\eqref{eq:assumption}. Does
  statement~\eqref{eq:assumption-new} hold in this case?
\end{problem}

\begin{problem}
  Let $\setA\subset\setM(N,M)$ and $\setB\subset\setM(M,N)$ be finite sets of
  matrices that satisfy assumption~\eqref{eq:tailineq}. Does
  statement~\eqref{eq:tailineq-new} hold in this case?
\end{problem}

\providecommand{\KeyWords}[1]{#1}

\end{document}